\documentclass[reqno,9 pt]{amsart}
\usepackage{amsmath,amsxtra,amssymb,latexsym, amscd,amsthm}
\usepackage[unicode]{hyperref}
\usepackage{array,tabularx,longtable,multicol,indentfirst,fancybox,color}
\usepackage{graphicx}
\usepackage{multicol}
\usepackage{mathrsfs}
\usepackage{enumerate}
\usepackage{tikz}
\usepackage{tikz-cd} 
\usepackage{tikz,tkz-tab}
\usetikzlibrary{calc, angles, quotes, intersections, positioning, patterns}
 
\advance\hoffset-1truecm\relax
\newtheorem{thm}{Theorem}[section]

\theoremstyle{definition}

\theoremstyle{remark}
\newtheorem{rem}{Remark}[section]
\numberwithin{equation}{section}
\numberwithin{equation}{section}
\begin{document}
	\title[Analog version of Hausdorff--Young's theorem for Quadratic Fourier transforms and boundedness of oscillatory integral]{Analog version of Hausdorff--Young's theorem for Quadratic Fourier\\ transforms and boundedness of oscillatory integral operator}
	
	\date{\today}%\today		
	\author[Trinh Tuan \& Lai Tien Minh]{Trinh Tuan$^{1*}$  \& Lai Tien Minh$^2$}
	\thanks{$^*$Corresponding author}

	\maketitle	
\begin{center}
	%\scriptsize
	$^1$Department of Mathematics, Faculty of Sciences, Electric Power University,\\ 235-Hoang Quoc Viet Rd., Bac Tu Liem Dist., Hanoi, Vietnam.\\
	E-mail: \texttt{tuantrinhpsac@yahoo.com}
	\\$^2$Department of Mathematics, Hanoi Architectural University,\\ Km. 10-Tran Phu Rd., Ha Dong Dist., Hanoi, Vietnam.\\
	E-mail: \texttt{minhlt@hau.edu.vn}
\end{center}	
\begin{abstract}
The purpose of this paper is twofold. The first aim is based on Riesz--Thorin's interpolation theorem, we prove new Hausdorff--Young type inequalities for the Quadratic Fourier transforms in (\textit{Ann. Funct. Anal.} 2014;5(1):10--23) and linear canonical transforms in (\textit{Mediterr. J. Math.} 2018;15,13), which were introduced by Castro et al. The second aim is to investigate the boundedness of the oscillatory integral operator with polynomial phases, which is also presented in the last section of the article. 
	\vskip 0.3cm
	\noindent\textbf{Keywords.} Hausdorff–Young theorem; Riesz--Thorin interpolation theorem; Quadratic Fourier transform; linear canonical transform; oscillatory integrals operator.
	\vskip 0.3cm
	
	\noindent \textbf{AMS Classifications.}  42A38; 42B10; 44A05; 26D10. 
\end{abstract}	
%\tableofcontents

\section{Introduction}
The classical Hausdorff--Young inequality for the Fourier transform \cite{3,1} is stated as follows: \textquotedblleft Suppose that $1\leq p\leq 2$, and $p_1$ is the conjugate exponent of $p$. If $f\in L_p(\mathbb{R}^n)$ then $Ff\in L_{p_1}(\mathbb{R}^n)$ and inequality $\|Ff\|_{L_{p_1}(\mathbb{R}^n)} \leq C\|f\|_{L_p(\mathbb{R}^n)}$ holds.$\textquotedblright$ This points out that the image of the Fourier transform is bounded by the original function on the conjugate exponential space. To study the next problems, the coefficient $C$ often chooses equal $1$. However, an outstanding result of Beckner  \cite{2}, following a special case of Babenko \cite{10} (when $p$ is an even number), showed that if $p$ is an arbitrary number in the interval [1,2], then: $\|Ff\|_{L_{p_1}} \leq \big[ \frac{p^{1\textfractionsolidus p}}{p_1^{1\textfractionsolidus p_1}}\big]^{\frac{n}{2}} \|f\|_{L_p}$, $\forall f\in L_p (\mathbb R^n)$. This is an improvement of the standard Hausdorff--Young inequality, as the context $p\leq 2$ and $p_1\geq2$ ensures that the number appearing on the right-hand side of this inequality is less than or equal to $1$. Moreover, this number cannot be replaced by a smaller one, since equality is achieved in the case of the Gaussian function $f: \mathbb{R} \rightarrow \mathbb{C}$ of the form $f(x) = \exp({-\pi|x|^2})$. In this sense, Beckner's result gives an optimal  \textquotedblleft sharp$\textquotedblright$ version of the Hausdorff--Young inequality \cite{2}. In the language of normed vector spaces, it says that the operator norm of the bounded linear mapping $L_p (\mathbb R^n)\rightarrow L_{p_1} (\mathbb R^n)$, as defined by the Fourier transform, is exactly equal to the constant $C=\big[ \frac{p^{1\textfractionsolidus p}}{p_1^{1\textfractionsolidus p_1}}\big]^{\frac{n}{2}}$.
Before presenting the main results of the article, we will recall some ideas about generalizing the Fourier transform previously proposed in \cite{5,6,4,11}.

First of all,  according to \cite{4}, an oscillating integral is introduced and studied by Phong-Stein, which has the forms
$(\mathcal{T}_\lambda\phi)(x) := \int_{\mathbb{R}} e^{i\lambda S(x,y)} \chi(y)\phi(y) dy,$ with $\phi \in C_0^{\infty}(\mathbb{R}),
$
where 
\begin{equation}\label{dfS}
S(x,y)=\displaystyle\sum_{k=1}^{n-1}\displaystyle\alpha_{k-1} x^{n-k} y^k
\end{equation}
is a homogeneous polynomial of degree $n$. The sharp decay estimate of the operator $(\mathcal{T}_\lambda)$ in the space $L_2(\mathbb{R})$ has been shown to be stable \cite{4}. The general oscillatory integral theory originates at the heart of harmonic analysis, in which Fourier’s case is the original and probably the best example of an oscillatory integral. It leads us to consider more general oscillatory integrals. There have been many efforts to estimate norm decay rates of Fourier oscillatory integrals \cite{6,phong1994models,stein1993harmonic}. However, it is difficult to generalize this result to general, higher-dimensional cases.

Secondly, in \cite{5}, based on the idea of a problem modeled by the heat equation by applying the theory of reproducing kernels \cite{11}, Castro et al. introduced and studied a Quadratic Fourier operator, which is considered global Quadratic functions in the exponent of the associated transforms. By choosing appropriate parameters in operator $(\mathcal{T}_\lambda)$, the authors introduced the following transform:
\begin{equation}\label{eq1.2}
f(x):= \int_{\mathbb{R}} e^{-ix(a\xi^2+b\xi +c)} \mathcal{F}(\xi) d\xi.
\end{equation}
In addition, by analyzing six subcases of \eqref{eq1.2} within a reproducing kernel Hilbert space framework gives us a comprehensive view of the Fourier transform with this type of phase function.
Later on, continuing the idea of generalizing the classical Fourier transform, in \cite{6} the authors considered the phase function, which has the form 
\begin{equation}\label{Q}
Q(x,y)=ax^2+bxy +cy^2+dx+ey\,\, ( a,b,c,d,e\in \mathbb{R},\,b\ne 0).
\end{equation}
The corresponding Quadratic Fourier operator $(\mathbb{Q})$ is given by \cite{6}:
$
(\mathbb{Q}f)(x):=\frac{1}{\sqrt{2\pi}}\int_{\mathbb{R}} e^{iQ(x,y)} f(y) dy.
$
Moreover, some basic properties for the integral operator $(\mathbb{Q})$ as well as some new convolutions associated with $(\mathbb{Q})$ will be introduced, which give us some interesting applications for the solvability of some classes of convolutional integral equations \cite{6}.

Being directly influenced by the above is also the driving motivation for this work. The first aim of this paper is to establish the Hausdorff--Young type theorems for several classes of Quadratic Fourier operators, introduced by Castro and co-workers in \cite{5,6}. Our achieved results are given by Theorem~\ref{thm2.1}, Theorem~\ref{thm2.2}, Theorem~\ref{thm2.2a}, and Theorem~\ref{thm2.3}, respectively. These obtained results lead to the boundedness of the Quadratic Fourier transforms in the conjugate space $L_{p_1}(\mathbb{R}),$ with $1\textfractionsolidus p + 1\textfractionsolidus p_1 =1,$ and $p \in [1,2]$. We would like to emphasize that the upper bound coefficients of all these inequalities are expressed specifically, but they are not the most optimally sharp version. In the last section, the boundedness of the oscillating integral operator $(\mathcal{T}_{\lambda})$ in the space $L_q(\mathbb{R})$ will be investigated, which is also the remaining goal of this article. Our obtained results are a continuation and improvement of the previous \cite{5,6} results. We also restate Minkowski's integral inequality as a useful tool to prove the main results later.
\begin{itemize}
	\item {\it Minkowski's integral inequality:} Follow \cite{12}, assume two measure spaces $(\theta_1,\mu_1)$ and $(\theta_2,\mu_2)$, and a measurable function
	$F(\cdot,\cdot): \theta_1\times\theta_2\rightarrow \mathbb{C}$, then
	\begin{equation}\label{eq1.3}
	\left[\int_{\theta_1}\left|\int_{\theta_1} F(x,y) d\mu_1(x)\right|^{s} d\mu_2(y)\right]^{\frac{1}{s}} \leq \int_{\theta_1}\left(\int_{\theta_2} |F(x,y)|^{s}d\mu_2(y)\right)^{\frac{1}{s}} d\mu_1(x),
	\end{equation}
	for any $s\geq 1$.
	
	\item According to (Theorem 12, p. 25 in \cite{9}): The formula
	\begin{equation}\label{eq1.4}
	\frac{1}{2}\left[f(x+0)+f(x-0)\right] = \lim_{\lambda\rightarrow \infty} \frac{1}{\pi}\int_{\mathbb{R}} f(t) \frac{\sin\lambda(x-t)}{x-t} dt
	\end{equation}
	holds true for if function $\frac{f(x)}{1+|x|}$ belongs to $L_1 (\mathbb{R})$,
	or
	$\frac{f(x)}{x}$ is of bounded variation in interval $(a, \infty)$ and $(-\infty,-a)$ for some positive $a$, and tends to $0$ at infinity.
\end{itemize}
%2
\section{Theorems of the Hausdorff--Young type}\label{sec2}
In this section, we will construct Hausdorff--Young type inequalities for the operators $(\mathbb{T}_\lambda)$, $(\mathbb{F}_k)$, $(\mathbb{H}_Q)$, which are defined by the formulas \eqref{eq2.1}, \eqref{eq2.3},  \eqref{eq2.3b},  and \eqref{eq2.5}, respectively. To achieve the above goals, we will use the techniques in \cite{5,6,4} together with the techniques of the Riesz--Thorin interpolation theorem (see Theorem~1.9,  p. 16 in \cite{3}).
%2.1
\subsection{Hausdorff--Young type theorem for the operator $(\mathbb{T}_\lambda)$}\label{sec2.1}
We will now consider the Hausdorff--Young type theorem for the operator $(\mathbb{T}_\lambda)$, which defined as follows \cite{6}:
\begin{equation}\label{eq2.1}
(\mathbb{T}_\lambda f)(x) = \int_{\mathbb{R}} e^{i\lambda Q(x,y)}\psi(x,y) f(y) dy,\quad f\in C_0^\infty(\mathbb{R}),
\end{equation}
where $\psi(x,y)$ is smooth compactly supported function on $\mathbb{R}^2$ and $Q(x,y)$ is taken by \eqref{Q}.
% by $Q(x,y):=ax^2+bxy+cy^2+dx+ey$ $(a,b,c,d,e\in\mathbb{R}, b\ne 0)$.
%thm2.1
\begin{thm}[Hausdorff--Young type inequality for $(\mathbb{T}_\lambda)$ transform]\label{thm2.1}
	Let $1\leq p\leq 2$ and $p_1$ be the conjugate exponential of $p$.  If $f\in L_p(\mathbb{R})$ then $(\mathbb{T}_\lambda f)\in L_{p_1}(\mathbb{R})$ and we have the following estimation
	\begin{equation}\label{eq2.2}
	\|\mathbb{T}_{\lambda}f\|_{L_{p_1}(\mathbb{R})} \leq \mathcal{C}_1\|f\|_{L_p(\mathbb{R})},\quad \mathcal{C}_1=C_1^{\frac{2}{p}-1}\left(\frac{C_1}{\lambda}\right)^{2\left(\frac{1}{p}-1\right)},
	\end{equation}
	where the constant $C_1$ is independent on $\lambda$.
\end{thm}
\begin{proof}
	Let $(Fg)(x)=\frac{1}{\sqrt{2\pi}}\int_{\mathbb{R}} e^{ixy} g(y) dy$ be the well-known Fourier transform.  Suppose that $M \subset \mathbb{R}^2$ is the compact support of $\psi(x,y)$.  Taking profit that $\psi(x,y)$ is uniformly bounded on $M\times M$. Thus,  there exists a constant $C_1$ such that
	$
	|\psi(x,y)| \leq C_1,\ \forall (x,y)\in M\times M.
	$
	For shortness of notation, let us consider $g_{\lambda}(y)=e^{i\lambda(cy^2+ey)}f(y)$.  A first direct computation yields $|g_{\lambda}(y)|=|f(y)|$.  Without loss of generality, 
	assume that $f\in L_p(\mathbb{R})$ and $\chi_M(x), \chi_M(y)$ stand for the characteristic functions with variable $x$ and $y$,  respectively.  We may observe that
	\begin{align*}
	\int_{\mathbb{R}} |\chi_M(y) \psi(x,y) g_{\lambda}(y)|^p dy &= \int_{\lambda} |\chi_M(y)\psi(x,y)|^p |g_{\lambda}(y)|^p dy\\
	& \leq C_1^p\int_{\mathbb{R}} |g_{\lambda}(y)|^p dy = C_1^p\|f\|_p^p< \infty.
	\end{align*}
	It follows that $\psi(x,y) g_{\lambda}(y)\in L_p(\mathbb{R})$.
	A direct computation give us
	\begin{align*}
	(\mathbb{T}_{\lambda}f)(x)&=e^{i\lambda(ax^2+dx)}\int_{\mathbb{R}} e^{ib\lambda xy}g_{\lambda}(y)\psi(x,y) dy
	= \sqrt{2\pi}e^{i\lambda(ax^2+dx)}F\{g_{\lambda}(y)\psi(x,y)\}(b\lambda x).
	\end{align*}
	Using Hausdorff--Young theorem for the Fourier transform $(F)$ (see \cite{2,3}),  we obtain
	$
	F\{g_{\lambda}(y)\psi(x,y)\}(b\lambda x)$ belongs to $L_{p_1}(\mathbb{R}).
	$
	Since $\big|e^{i\lambda(ax^2+dx)}\big|=1$,  this indicates that $(\mathbb{T}_{\lambda}f)\in L_{p_1}(\mathbb{R})$. On the other hands
	\begin{align*}
	|(\mathbb{T}_{\lambda}f)(x)| &\leq \int_{\mathbb{R}} \left|e^{i\lambda Q(x,y)}\right| |f(y)| |\psi(x,y)| \chi_M(y) dy\\
	&= \int_{\mathbb{R}} |f(y)| |\psi(x,y)| \chi_M(y) dy\\
	&\leq C_1\int_{\mathbb{R}} |f(y)| dy
	=C_1\|f\|_{L_1(\mathbb{R})}<\infty\,,\forall f \in L_1(\mathbb{R}).
	\end{align*}
	%Thus,  if $f\in L_1(\mathbb{R})$ then $|(\mathbb{T}_{\lambda}f)(x)| \leq C_1\|f\|_1$.  
	Therefore,  the estimation above can be recast as
	$
	\|\mathbb{T}_{\lambda}f\|_{L_1(\mathbb{R})} := \underset{x\in M}{\mathrm{ess\,sup}}|(\mathbb{T}_{\lambda}f)(x)| \leq C_1\|f\|_{L_1(\mathbb{R})}
	$ is finite.
	Therefore,  we deduce that $(\mathbb{T}_{\lambda}f)$ is bounded operator from $L_1(\mathbb{R})\rightarrow L_{\infty}(\mathbb{R})$.  Moreover, Theorem~4.4 in \cite{6} gives us
	$
	\|\mathbb{T}_{\lambda}\|_{L_2(\mathbb{R})}\leq \frac{C_1}{\sqrt{\lambda}}.
	$
	Thus, the operator $(\mathbb{T}_{\lambda})$ is bounded from $L_2(\mathbb{R})\rightarrow L_2(\mathbb{R})$.  Applying the Riesz--Thorin interpolation theorem (see \cite{3}, Theorem~1.19, p. 16) we receive that the operator $(\mathbb{T}_{\lambda})$ is bounded from $L_p( \mathbb{R})\rightarrow L_{p_1}(\mathbb{R})$, where $1\leq p\leq 2$, $\frac{1}{p}+\frac{1}{p_1} =1$ and $\frac{1}{p}=\frac{\alpha}{1}+\frac{1-\alpha}{2}$ (and then $\alpha=\frac{2}{ p}-1$). Moreover,  we have the following estimation
	$
	\|\mathbb{T}_{\lambda}f\|_{L_{p_1}(\mathbb{R})} \leq C_1^{\alpha}\left(\frac{C_1}{\sqrt{\lambda}}\right)^{1-\alpha} \|f\|_{L_p(\mathbb{R})},\ (0<\alpha<1),
	$
	which is the desired result.
	%This is equivalent to
	%$$
	%\|\mathbb{T}_{\lambda}f\|_{L_{p_1}(\mathbb{R})} \leq C_1^{\alpha}\left(\frac{C_1}{\sqrt{\lambda}}\right)^{1-\alpha}\|f\|_{L_p(\mathbb{R})}.
	%$$
\end{proof}
We now shall retrieve the some special cases of the norm inequality  \eqref{eq2.2} as follows
%rem2.1
\begin{rem}\label{rem}
	When $p=1$, then $\|\mathbb{T}_{\lambda} f\|_{L_{\infty}(\mathbb{R})} \leq C_1\|f\|_{L_1(\mathbb{R})},\,\, \forall f\in L_1(\mathbb{R}).$
	Let $p=2$.  In such aparticular case we realize that $\|\mathbb{T}_{\lambda}f\|_{L_2(\mathbb{R})} \leq \frac{C_1}{\sqrt{\lambda}} \|f\|_{L_2(\mathbb{R})},$ which is therefore the estimation in \cite{6} (Theorem~4.4, p. 13).
\end{rem}
%2.2
\subsection{Hausdorff--Young type theorems for the operator $(\mathbb{F}_k)$}\label{sec2.2}
Inspired by the operator \eqref{eq1.2} in \cite{5}.  More globally,  let us consider the operators $(\mathbb{F}_k)$ with $k\in\{1,2\}$,  which are defined as follows
\begin{equation}\label{eq2.3}
(\mathbb{F}_1f)(x) = \int_{\mathbb{R}} e^{-ix(ay^2+by+c)} \psi(x,y) f(y) dy,
\end{equation}
\begin{equation}\label{eq2.3b}
(\mathbb{F}_2f)(x) = \int_{\mathbb{R}} e^{-i(ax^2+bx+c)y} \psi(x,y) f(y) dy,
\end{equation}
where $\psi(x,y)$ is smooth compactly supported function on $\mathbb{R}^2$.  We have the following results
%thm2.2
\begin{thm}[Hausdorff--Young type inequality for $(\mathbb{F}_1)$ transform]\label{thm2.2}
	If $p\in [1,2]$ and $\frac{1}{p}+\frac{1}{p_1}=1$.  Suppose that $f\in L_p(\mathbb{R})$ then $(\mathbb{F}_1f) \in L_{p_1}(\mathbb{R})$. Moreover, the following estimation holds 
	\begin{equation}\label{eq2.4}
	\|\mathbb{F}_1f\|_{L_{p_1}(\mathbb{R})} \leq \mathcal{C}_2\|f\|_{L_p(\mathbb{R})},
	\end{equation}
	where $\mathcal{C}_2=C_1R^{2\left(\frac{1}{p}-1\right)}$ and $R$ is the diameter of $M$, where $M$ is the compact support of function $\psi(x,y)$.
\end{thm}

\begin{proof}
	We infer directly that
	$
	(\mathbb{T}f)(x) = e^{-icx}\int_{\mathbb{R}}e^{-ibxy}\left[e^{-iaxy^2}f(y)\psi(x,y)\right] dy.
	$
	By using a similar way to the proof of Theorem~\ref{thm2.1},  we achieve
	\begin{align*}
	\int_{\mathbb{R}}\left|\chi_M(y) e^{-iaxy^2} f(y)\psi(x,y)\right|^p dy &= \int_{\mathbb{R}}\left|\chi_M(y)\psi(x,y)\right|^p |f(y)|^p dy \\
	&\leq C_1^p\int_{\mathbb{R}}|f(y)|^p dy.
	\end{align*}
	From the assumption $f\in L_p(\mathbb{R})$,  it is easy to check that
	$
	\int_{\mathbb{R}} \left|\chi_M(y)e^{-iaxy^2}f(y) \psi(x,y)\right|^p dy \leq C_1^p\|f\|_{p}^p 
	$ is finite, 
	this leads $e^{-iaxy^2}f(y)\psi(x,y)\in L_p(\mathbb{R})$.  Obviously,  we can write
	$$
	(\mathbb{F}_1f)(x) = \sqrt{2\pi} e^{-icx} F\left\{e^{-iaxy^2}f(y)\psi(x,y)\right\}(-bx).
	$$
	Using the Hausdorff--Young theorem for the Fourier transform $(F)$ (see \cite{2,3}),  we arrive at $$F\left\{e^{-iaxy^2}f(y)\psi(x,y)\right\}(bx)\in L_{p_1}(\mathbb{R})$$ with $(1\textfractionsolidus p + 1\textfractionsolidus p_1 =1,$ and $1\leq p\leq 2)$.  Hence,  this implies that $(\mathbb{F}_1f)(x) \in L_{p_1}(\mathbb{R})$.
	By performing a similar way as Theorem~\ref{thm2.1},  we conclude that $(\mathbb{F}_1f)$ is bounded operator from $L_1(\mathbb{R})\rightarrow L_{\infty}(\mathbb{R})$. Moreover, we can get
	$
	\|\mathbb{F}_1f\|_{L_{\infty}(\mathbb{R})} \leq C_1\|f\|_{L_1(\mathbb{R})}.
	$
	Using Minkowski's integral inequality \eqref{eq1.3} and the Cauchy-Bunyakovsky-Schwarz inequality, we obtain
	\begin{align*}
	\int_{\mathbb{R}} |(\mathbb{F}_1f)(x)|^2 dx&= \int_{\mathbb{R}}\left|\int_{\mathbb{R}} e^{-ix(ay^2+by+c)} \psi(x,y) f(y) dy\right|^2 dx\\
	&\leq \left[\int_{\mathbb{R}}\left(\int_{\mathbb{R}}\left|\chi_M(x)\chi_M(y)e^{-ix(ay^2+by+c)} \psi(x,y) f(y)\right|^2dx\right)^{\frac{1}{2}}dy\right]^{2}\\
	&=\left(\int_{\mathbb{R}}|f(y)| \chi_M(y) dy\right)^2\left(\int_{\mathbb{R}}\chi_M(x)|\psi(x,y)|^2 dx\right)\\
	&\leq \left(\int_M\chi_M(y)dy\right)\left(\int_M|f(y)|^2 dy\right)\left(\int_{\mathbb{R}} \chi_M(x)|\psi(x,y)|^2 dx\right)\\
	&\leq C_1^2\left(\int_M\chi_M(y) dy\right)\left(\int_M \chi_M(x) dx\right)\left(\int_{\mathbb{R}} |f(y)|^2 dy\right)
	\leq C_1^2R^2\|f\|_{L_2(\mathbb{R})}^2,
	\end{align*}
	where $R=\displaystyle\max_{(x,y)\in M} \rho(x,y)$.  Therefore
	$
	\|\mathbb{F}_1f\|_{L_2(\mathbb{R})} \leq C_1R\|f\|_{L_2(\mathbb{R})}$ is finite, $\forall f\in L_2(\mathbb{R}).
	$
	Thus, $(\mathbb{F}_1)$ is the bounded operator from $L_2(\mathbb{R})\rightarrow L_2(\mathbb{R})$.
	Applying the Riesz–Thorin interpolation theorem \cite{4}, we deduced that the operator $(\mathbb{F}_1)$ is bounded from $L_p(\mathbb{R})\rightarrow L_{p_1}(\mathbb{ R})$ with $\frac{1}{p}+\frac{1}{q}=1$, and $p \in [1,2]$.  Moreover, since $\frac{1}{p}=\frac{\alpha}{1}+\frac{1-\alpha}{2},$ infer that $\alpha=\frac{2}{p}-1$. Then,  we have the desired decay rate
	$\|\mathbb{F}_1f\|_{L_{p_1}(\mathbb{R})} \leq C_1^{\alpha}(C_1R)^{1-\alpha}\|f\|_{L_p(\mathbb{R})},$
	which yields the desired conclusion of inequality. The proof is completed.
\end{proof}
The following theorem will be omitted because its proof is very similar to the proof of Theorem \ref{thm2.2}.
%Similarly, using formula \eqref{eq2.3b},  we obtain
\begin{thm}[Hausdorff--Young type inequality for $(\mathbb{F}_2)$ transform]\label{thm2.2a}
	Suppose that $1\leq p\leq 2$ satisfy condition $\frac{1}{p}+\frac{1}{p_1}=1$. For any functions $f\in L_p(\mathbb{R})$.  Then $(\mathbb{F}_2f) \in L_{p_1}(\mathbb{R})$. Furthermore,  the following inequality holds
	\begin{equation}\label{eq2.4c}
	\|\mathbb{F}_2f\|_{L_{p_1}(\mathbb{R})} \leq \mathcal{C}_2\|f\|_{L_p(\mathbb{R})},
	\end{equation}
	where $\mathcal{C}_2=C_1R^{2\left(\frac{1}{p}-1\right)}.$
	%where $R$ is the diameter of $M$ - compact price of $\psi(x,y)$.
\end{thm}
%2.3
\subsection{Hausdorff--Young type theorem for the operator $(\mathbb{H}_Q)$}\label{sec2.3}
In this section, we construct the Hausdorff--Young type theorem for the Quadratic Hartley operator $(\mathbb{H}_Q)$, which is defined as follows
\begin{equation}\label{eq2.5}
(\mathbb{H}_Qf)(x)=\sqrt{\frac{|b|}{2\pi}} \int_{\mathbb{R}} \left[\frac{e^{-iQ(x,y)} + ze^{-iQ(-x,y)}}{2}\right]f(y) dy,\, b\ne 0,~z\in\left\{\sqrt{-1}\right\},
\end{equation}
where $Q(x,y)$ is polynomial taken by \eqref{Q}. 
In the following, some special cases of $(\mathbb{H}_Q)$ transform will presented in the remark
\begin{rem}
	Let $d=e=0$.  In such a particular case,  $(\mathbb{H}_Q)$ transform becomes the canonical Hartley transform (refer \cite{7}).
	If $a=b=c=d=e=0$, $b=1$ and,  $z=i$ then we would like to notice that $(\mathbb{H}_Q)$ is simply the well-knows Hartley transform \cite{8,MedJmath2022}.
\end{rem}
%thm2.3
\begin{thm}[Hausdorff--Young-type inequality for the  $(\mathbb{H}_Q)$ transform]\label{thm2.3}
	Let $p \in [1,2]$ and $p_1$ be the conjugate exponent of $p$. Then, we have the following estimation
	\begin{equation}\label{eq2.6}
	\|\mathbb{H}_Qf\|_{L_{p_1}(\mathbb{R})} \leq \mathcal{C}_3 \|f\|_{L_p(\mathbb{R})},\quad \forall f\in L_p(\mathbb{R}),
	\end{equation}
	where $\mathcal{C}_3=\left(\sqrt{\frac{|b|}{2\pi}}\right)^{\alpha}$, $\alpha=\frac{2}{p}-1$.
\end{thm}
\begin{proof}
	Without loss of generality,  we assume that $z=i$.  Since $g_1(y)=e^{-i(cy^2+ey)} f(y)$ then $|g_1(y)|=|f(y)|$. By simple computations, we have
	\begin{equation}\label{H1}
	\begin{aligned}
	(\mathbb{H}_Qf)(x)&= \frac{\sqrt{|b|}}{2} e^{-i(ax^2+dx)} \frac{1}{\sqrt{2\pi}}\int_{\mathbb{R}} e^{-ibxy} g_1(y) dy
	+ \frac{\sqrt{|b|}}{2} e^{-i(ax^2-dx)} \frac{1}{\sqrt{2\pi}} \int_{\mathbb{R}} e^{ibxy} g_1(y) dy\\
	&=\frac{\sqrt{|b|}}{2} e^{-i(ax^2+dx)}(Fg_1)(-bx) + \frac{\sqrt{|b|}}{2}e^{-i(ax^2-dx)}(Fg_1)(bx).
	\end{aligned}
	\end{equation}
	Because of $f\in L_p(\mathbb{R})$,  it is straightforward to get $g_1\in L_p(\mathbb{R})$ $(1\leq p\leq 2)$.  Moreover, making use of the Hausdorff--Young inequality (\cite{2,3}), we obtain  
	$(Fg_1)(bx) \in L_{p_1}(\mathbb{R})$, and  $(Fg_1)(-bx) \in L_{p_1}(\mathbb{R}), $
	where $\frac{1}{p}+\frac{1}{p_1}=1,$ and  $p\in [1,2]$. 
	It is easy to verify that $\big|e^{-i(ax^2+dx)}\big| = \big|e^{-i(ax^2-dx)}\big|=1$.  Now,  owing to \eqref{H1}, we derive $(\mathbb{H}_Qf)\in L_{p_1}(\mathbb{R}).$
	On the other hand, we proceed the kernel of $(\mathbb{H}_Q)$ as
	\begin{align*}
	\frac{e^{-iQ(x,y)}+ie^{-iQ(-x,y)}}{2}&= \frac{e^{-iQ(x,y)} + \frac{1+i}{1-i}e^{-iQ(-x,y)}}{2}\\
	&= \frac{e^{-i(ax^2+cy^2+ey)}}{1-i}\left[\frac{e^{-i(by+d)x}+e^{i(by+d)x}}{2}+ i\frac{e^{i(by+d)x}-e^{-i(by+d)x}}{2}\right]\\
	&=\frac{e^{-i(ax^2+cy^2+ey)}}{1-i}\big[\cos(by+d)x - \sin(by+d)x\big].
	\end{align*}
	Therefore
	$
	\left|\frac{e^{-iQ(x,y)}+ie^{-iQ(-x,y)}}{2}\right| = \frac{1}{\sqrt{2}} \left|\cos(by+d)x-\sin(by+d)x\right|\leq 1,\ \forall x,y.
	$
	Then,  using the above relation, we obtain
	\begin{align*}
	|(\mathbb{H}_Qf)(x)| &\leq \sqrt{\frac{|b|}{2\pi}} \int_{\mathbb{R}}\left|\frac{e^{-iQ(x,y)}+ze^{-iQ(-x,y)}}{2}\right||f(y)| dy\\
	&\leq \sqrt{\frac{|b|}{2\pi}} \int_{\mathbb{R}}|f(y)|dy = \sqrt{\frac{|b|}{2\pi}}\|f\|_{L_1(\mathbb{R})} < \infty,\quad \forall f\in L_1(\mathbb{R}).
	\end{align*}
	Moreover
	$
	\|\mathbb{H}_Qf\|_{L_\infty(\mathbb{R})} := \underset{x\in\mathbb{R}}{\mathrm{ess\,sup}}|(\mathbb{H}_Qf)(x)| \leq \sqrt{\frac{|b|}{2\pi}} \|f\|_{L_1(\mathbb{R})}.
	$
	Thus, we infer that the $(\mathbb{H}_Q)$ is bounded operator from $L_1(\mathbb{R})\rightarrow L_\infty(\mathbb{R})$.
	Next,  we turn to  prove that $(\mathbb{H}_Q)$ is the unitary operator in $L^2(\mathbb{R}).$ Suppose that $b>0$, we then have
	\begin{align}
	\|\mathbb{H}_Qf\|_{L_2(\mathbb{R})}^2&=\langle \mathbb{H}_Qf,\mathbb{H}_Qf\rangle = \int_{\mathbb{R}} (\mathbb{H}_Qf)(x) \overline{(\mathbb{H}_Qf)(x)} dx\notag\\
	&=\frac{b}{4\pi}\! \int_{\mathbb{R}^3}\left[e^{-iQ(x,y)}+z e^{-iQ(-x,y)}\right]\!\!\left[e^{iQ(x,y)}+\bar{z}e^{iQ(-x,y)}\right]\!\!f(y)\bar{f}(u) dy du dx\notag\\
	&=\frac{b}{4\pi}\int_{\mathbb{R}^3}\left\{\left[e^{-[Q(x,y)-Q(x,y)]} + e^{i[Q(-x,y)-Q(-x,y)]}\right]\right.\notag\\
	&\qquad\qquad\left.+z\left[e^{-i[Q(x,y)-Q(-x,u)]} - e^{-i[Q(-x,y)-Q(x,u)]}\right]\right\} f(y)\bar{f}(u) dy du dx\notag\\
	&=\frac{b}{2\pi}\int_{\mathbb{R}^2} \left\{\int_{\mathbb{R}}\left[\frac{e^{-ibx(y-u)} + e^{ibx(y-u)}}{2}\right.\right.\notag\\
	&\hspace{1.5cm}\left.\left. + z\frac{e^{-i[b(y+u)+2d]x} - e^{i[b(y+u)+2d]x}}{2}\right] dx\right\} g_1(y)\overline{g_1}(u) dy du\notag\\
	&=\frac{b}{2\pi}\int_{\mathbb{R}^2} g_1(y)\overline{g_1}(u) \left\{\int_{\mathbb{R}}\left(\cos b(y-u)x -z\sin[b(y+u)+2d]x\right)dx\right\} dy du\notag\\
	&=\frac{b}{2\pi}\lim_{t\rightarrow\infty}\!\int_{\mathbb{R}^2}g_1(y) \overline{g_1}(u)\!\left\{\!\int_{-t}^t\!\!\! (\cos b(y-u)x - z\sin[b(y+u)+2d]x)dx\right\} dy du.\label{H3}
	\end{align}
	Thanks to formula $\displaystyle\int_{-t}^t \sin[b(y+u)+2d]x dx=0, $ the relation \eqref{H3} becomes
	$$
	\|\mathbb{H}_Qf\|_{L_2(\mathbb{R})}^2 = \int_{\mathbb{R}} g_1(y)\left\{\lim_{t\rightarrow\infty} \frac{1}{\pi}\int_{\mathbb{R}}\overline{g_1}(u) \frac{\sin bt(y-u)}{y-u} du\right\} dy.
	$$
	By virtue of the formula \eqref{eq1.4},  we derive
	$$
	\|\mathbb{H}_Qf\|_{L_2(\mathbb{R})}^2 = \int_{\mathbb{R}} g_1(y)\overline{g_1}(y) dy = \int_{\mathbb{R}} f(y)\bar{f}(y) dy = \int_{\mathbb{R}} |f(y)|^2dy = \|f\|_{L_2(\mathbb{R})}^2.
	$$
	Thus
	\begin{equation}\label{H2}
	\|\mathbb{H}_Qf\|_{L_2(\mathbb{R})} = \|f\|_{L_2(\mathbb{R})}<\infty, \forall f\in L_2(\mathbb{R}).
	\end{equation}
	In the case $b<0$,  we can be derived \eqref{H2} in the similar way.  Therefore, we obtain  $(\mathbb{H}_Q)$ transform is an isometric isomorphism (unitary) mapping in $L_2(\mathbb{R})$.  Hence,  $(\mathbb{H}_Q)$ is bounded operator from $L_2(\mathbb{R} )\rightarrow L_2(\mathbb{R})$.
	Applying the Riesz--Thorin interpolation theorem \cite{3},  we deduce that $(\mathbb{H}_Q)$ is the bounded operator from $L_p(\mathbb{R})\rightarrow L_{p_1}(\mathbb{R})$, where $p\in [1,2]$, $\frac{1}{p}+\frac{1}{p_1}=1$ and $\frac{1}{p}=\frac{\alpha}{1}+\frac{1-\alpha}{2}$ (and then $\alpha=\frac{2}{p}-1$). Therefore,  the desired estimation can be achieved as follows
	$
	\|\mathbb{H}_Qf\|_{L_{p_1}(\mathbb{R})} \leq \big(\sqrt{\frac{|b|}{2\pi}}\big)^{\frac{2}{p}-1} \|f\|_{L_p(\mathbb{R})}.
	$
	This indicates that \eqref{eq2.6} is proved.
\end{proof}

\begin{rem}
	If $p=1$, then estimate \eqref{eq2.6} becomes
	$
	\|\mathbb{H}_Qf\|_{L_{\infty}(\mathbb{R})} \leq \sqrt{\frac{|b|}{2\pi}} \|f\|_{L_1(\mathbb{R})}.
	$
	For case $p=2$,  we infer directly that
	$
	\|\mathbb{H}_Qf\|_{L_2(\mathbb{R})} = \|f\|_{L_2(\mathbb{R})}.
	$ Actually,  the coefficients $\mathcal{C}_1$, $\mathcal{C}_2$, and $\mathcal{C}_3$ in the estimates \eqref{eq2.2}, \eqref{eq2.4}, \eqref{eq2.4c},  and \eqref{eq2.6} are not the most optimally sharp. The problem of finding the best coefficients for the above estimates is still an open question.
\end{rem}
 
%3
\section{Boundedness of oscillatory integral operator $(\mathcal{T}_\lambda)$ in the $L_q(\mathbb{R})$ space}\label{sec3}
Let $\chi$ is a fixed $C^\infty$ cut-off function on $\mathbb{R}$ with compact supported.  Assume that  $S(x,y)$ is a homogeneous polynomial of degree $n$ in $(x,y)\in \mathbb{R}^2$,  which can be given by \eqref{dfS}.  The oscillatory integral with polynomial phases of $\phi$ is defined by \cite{4} as follows
\begin{equation}\label{eq3.1}
(\mathcal{T}_\lambda \phi)(x) = \int_{\mathbb{R}} e^{i\lambda S(x,y)} \chi(y) \phi(y) dy,\quad \phi\in C_0^\infty(\mathbb{R}).
\end{equation}
In the first place,  it is easy to see that \eqref{eq2.1}, \eqref{eq2.3}, \eqref{eq2.3b},  and \eqref{eq2.5} are not special cases of \eqref{eq3.1}. More importantly,  when the coefficients $\alpha_k$ of the polynomial $S(x,y)$ satisfy the given conditions,  in \cite{4} Phong and Stein introduced the estimation of the norm in space $L_{2}(\mathbb{R})$ as follows
$\|\mathcal{T}_\lambda\|_{L_2(\mathbb{R})} \leq C_1\lambda^{-\frac{1}{n}}.$
This section is devoted to proving the boundedness of the operator $(\mathcal{T}_\lambda \phi)$ in the space $L_{p_1}(\mathbb{R})$ when $\phi\in L_{p}(\mathbb{R})$, where $p_1$ is the conjugate exponent of $p$.  Our key result can be obtained by using Minkowski inequality \eqref{eq1.3} and H\"older inequality,  which is presented by 
%thm3.1
\begin{thm}
	Let $p, p_1$ be real numbers  in open interval $(1,\infty)$ such that $\frac{1}{p}+\frac{1}{p_1}=1$.  Then,  operator $(\mathcal{T}_\lambda \phi)$ extends to a bounded operator on $L_{p_1}(\mathbb{R})$ when $\phi\in L_{p}(\mathbb{R})$.
\end{thm}
\begin{proof}
	Without loss of generality, suppose that $M\subset \mathbb{R}^2$ is the compact support price of $\psi$ and $\chi_M(x)$, $\chi_M(y)$ is the characteristic function of the two variables $x$ and $y$. Then, 
	with the aid of the Minkowski inequality \eqref{eq1.3}, we can write
	\begin{align}
	\int_{\mathbb{R}} \left|\left(\mathcal{T}_\lambda \phi\right)(x)\right|^{p_1} dx&= \int_{\mathbb{R}} \left|\int_{\mathbb{R}} e^{i\lambda s(x,y)} \chi_M(x) \chi_M(y) \psi(x,y) \phi(y) dy\right|^{p_1} dx\notag\\
	&\leq \left[\int_{\mathbb{R}} \left(\int_{\mathbb{R}}\left|\chi_M(x) \chi_M(y) \psi(x,y) \phi(y)\right|^{p_1} dx\right)^{\frac{1}{p_1}} dy\right]^{p_1}\notag\\
	&= \left[\int_{\mathbb{R}}\left(\int_{\mathbb{R}} \chi_M(x) \chi_M(y) |\psi(x,y)|^{p_1} |\phi(y)|^{p_1} dx\right)^{\frac{1}{p_1}} dy\right]^{p_1}\notag\\
	&= \left[\int_{\mathbb{R}} |\phi(y)| \chi_M(y) dy\right]^{p_1} \left[\int_{\mathbb{R}}\chi_M(x) |\psi(x,y)|^{p_1} dx\right].\label{*}
	\end{align}
	Having now in mind that $|\psi(x,y)|\leq C_1$, $\forall (x,y)\in M\times M$, we realize that
	\begin{equation}\label{**}
	\int_{\mathbb{R}} \chi_M(x)|\psi(x,y)|^{p_1} dx\leq C_1^{p_1} <\infty.
	\end{equation}
	By combining \eqref{*} and \eqref{**},  we deduce that
	$
	\int_{\mathbb{R}} |(\mathcal{T}_\lambda\phi)(x)|^{p_1} dx \leq C_1^{p_1}\left[\int_{\mathbb{R}} |\phi(y)| \chi_M(y) dy\right]^{p_1}.
	$
	Using the H\"older inequality, with $1\textfractionsolidus p + 1\textfractionsolidus p_1 =1,$ and  $p,p_1\in (1,\infty)$, then we have
	$$
	\int_{\mathbb{R}}|(\mathcal{T}_\lambda\phi)(x)|^{p_1} dx \leq C_1^{p_1}\left[\left(\int_M |\phi(y)|^{p} dy\right)^{\frac{1}{p}} \left(\int_M \chi_M^{p_1}(y) dy\right)^{\frac{1}{p_1}}\right]^{p_1}.
	$$
	Since $\phi\in L_{p}(\mathbb{R})$,  we achieve
	$
	\|\mathcal{T}_\lambda\phi\|_{L_{p_1}(\mathbb{R})} \leq C_1\left(\int_{M} \chi_M^{p_1}(y)dy\right)^{\frac{1}{p_1}}\|\phi\|_{L_{p}(\mathbb{R})}
	$ is finite, 
	which prove $(\mathcal{T}_\lambda\phi)$ belongs to $L_{p_1}(\mathbb{R})$.  The proof is completed.
\end{proof}

\noindent Based on the analogous method as in Theorem \ref{thm2.1} when $p=1$ and $p_1=\infty$,  we also conclude that $(\mathcal{T}_\lambda)$ is bounded operator from $L_1(\mathbb{R} )\rightarrow L_{\infty}(\mathbb{R})$.  It would be interesting to obtain Hausdorff--Young theorem for the oscillating integral $(\mathcal{T}_\lambda)$.  However, $S(x,y)$ is a homogeneous polynomial of degree $n$ is one of the most challenging of proving this theorem.  Hence,  up to now,  constructing the Hausdorff--Young type theorem for the $(\mathcal{T}_\lambda)$ is still an open question.

\vskip 0.3cm
\noindent \textbf{Data availability} \\
Our manuscript has no associated data.

\vskip 0.3cm
\noindent \textbf{Funding}\\
\noindent This research received no specific grant from any funding agency.

\vskip 0.3cm
\noindent \textbf{Conflict of interest}\\
\noindent The authors have no conflicts of interest to declare that are relevant to the content of this article.
\vskip 0.3cm

\noindent \textbf{ORCID}\\
\noindent \textit{Trinh Tuan}  \url{https://orcid.org/0000-0002-0376-0238}\\
\noindent \textit{Lai Tien Minh} 
\url{https://orcid.org/0000-0003-2656-8246}

%\bibliographystyle{plain}
%\bibliography{ref_preprint_JOTH}
\end{document}